\documentclass[12pt]{article}
\usepackage{amssymb,latexsym,amscd,amsmath,amsfonts,enumerate,supertabular}
\usepackage{graphicx}
\usepackage{tabularx}
\usepackage{caption}
 \usepackage{epstopdf}
\newtheorem{Theorem}{Theorem}[section]

\newtheorem{Lemma}[Theorem]{Lemma}

\numberwithin{equation}{section}
\begin{document}

\title{{\bf Existence results and numerical solution for the Dirichlet problem for fully fourth order nonlinear equation
}}

%
%
\author{ Dang Quang A$^{\text a}$,  Nguyen Thanh Huong$^{\text b}$\\
$^{\text a}$ {\it\small Center for Informatics and Computing, Vietnam Academy of Science and Technology}\\
{\it\small (VAST),  18 Hoang Quoc Viet, Cau Giay, Hanoi, Vietnam}\\
{\small Email address: dangquanga@cic.vast.vn}\\
$^{\text b}$ {\it\small College of Sciences, 
Thainguyen University,
Thainguyen, Vietnam}\\
{\small Email address: nguyenthanhhuong2806@gmail.com}}
\date{ }
\maketitle


\begin{abstract}
\small  
{In this paper we study the existence and uniqueness of a solution and propose an iterative method for solving a beam problem which is described by the fully fourth order equation 
$$u^{(4)}(x)=f(x,u(x),u'(x),u''(x),u'''(x)), \quad 0 < x < 1$$ 
associated with the Dirichlet boundary conditions.\\
This problem was studied by several authors. Here we propose a novel approach by the reduction of the problem to an operator equation for the triplet of  the nonlinear term
$\varphi (x)=f(x,u(x),u'(x),u''(x),u'''(x))$ and the unknown values $u''(0), u''(1).$
Under some easily verified conditions on the function $f$ in a specified bounded domain, we prove the contraction of the operator. This guarantees the existence and uniqueness of a solution and the convergence of an iterative method for finding it. Many examples demonstrate the applicability of the theoretical results and the efficiency of the iterative method. The advantages of the obtained results over those of Agarwal are shown on some examples. }

\end{abstract}

\section {Introduction}
In the last decade the fully fourth order nonlinear differential equation
\begin{equation}\label{eqMain}
u^{(4)}(x)=f(x,u(x),u'(x),u''(x),u'''(x))
\end{equation}
has attracted a great attention from many researchers because it is the mathematical model of numerous engineering problems. Depending on the concrete engineering problems the equation \eqref{eqMain} is considered with different boundary conditions. Below we mention some works concerning typical boundary conditions.\par
First, it is worth to mention a very recent work of Y. Li \cite{Li2}. In this work the author consider the problem
\begin{equation}\label{BC0}
\begin{aligned}
u^{(4)}(x)&=f(x,u(x),u'(x),u''(x),u'''(x)), \quad 0 < x < 1, \\
u(0)&=u'(0)=u''(1)=u'''(1)=0,
\end{aligned}
\end{equation}
which models a statistically elastic beam fixed at the left and freed at the right end. Under several complicated conditions such as the growth conditions on infinity, a Nagumo-type condition, using the fixed point index theory in cones the author obtained some results of existence of positive solutions. Freeing the above conditions but requiring the Lipschitz condition in a specific bounded domain in \cite{AQ} we established the existence and uniqueness of a solution and proposed an iterative method for finding the solution. The method used is the reduction of the boundary value problem to an operator equation for the unknown function $\varphi (x)= f(x,u(x),u'(x),u''(x),u'''(x))$. This method is first proposed in our work \cite{ALQ}. A problem for the equation \eqref{eqMain} with somewhat different from those in \eqref{BC0} boundary conditions, namely, 
\begin{equation}\label{BC1}
u(0)=u''(0)=u'(1)=u'''(1)=0
\end{equation}
was studied in \cite{Min} by the lower and upper solution method and degree theory. In \cite{Bai} Bai also used the lower and upper solution method under the Nagumo  and monotonocity conditions proved the existence of solutions of the equation \eqref{eqMain} associated with boundary conditions \eqref{BC1}.
An another problem for the equation \eqref{eqMain} and the boundary conditions 
\begin{equation}\label{BC2}
u(0)=u''(0)=u(1)=u''(1)=0,
\end{equation}
which models deformation of an elastic beam with simply supported ends, was investigated in \cite{Li1}. There by the Fourier analysis method and the Leray-Schauder fixed point theorem the author established the existence of a solution under the linear growth condition. The uniqueness of the solution is guaranteed if adding the Lipschitz condition. In \cite{Dang3} we relaxed these conditions by the requirement only of the Lipschitz condition in a bounded domain. \par
Except for the boundary conditions in \eqref{BC0} and \eqref{BC1}, \eqref{BC2} in \cite{Feng}, \cite{Pei2} the authors considered the equation \eqref{eqMain} with the boundary conditions
\begin{equation}\label{BC3}
 u(0)=0, u'(1)=0, au''(0)-bu'''(0)=0, cu''(1)+du'''(1)=0. 
 \end{equation} 
The existence of a solution was established by the Leray-Schauder degree theory and by means of the lower and upper solution method. A constructive proof of the existence of a solution also was given in \cite{Du}.\par 

As was seen above, there are many works devoted to boundary value problems for the fully fourth order nonlinear equations, but the boundary conditions associated are not the Dirichlet ones.
To the best of our knowledge, there is only a work \cite{Agarwal} concerning the existence for the problem with Dirichlet boundary conditions. In this work the authors considered the Dirichlet problem
\begin{equation}\label{eq0}
\begin{aligned}
u^{(4)}(x)=f(x,u(x),u'(x),u''(x),u'''(x)), \quad a < x < b, \\
u(a)=A_1, \quad u(b)=B_1, \quad u'(a)=A_2, \quad u'(b)=B_2,
\end{aligned}
\end{equation}
where $f \in C[[a,b]\times \mathbb{R}^4\rightarrow \mathbb{R}]$, $A_i$, $B_i$ $(i=1,2)$ are real constants. Some theorems on existence of a solution based on the use of Shauder fixed point theorem were established and the convergence of Picard iterations were proved. But as the examples in this paper showed that not always these theorems are applicable, and for the convergence of Picard iterations it requires the Lipschitz condition posed on the nonlinear function $f$ in a domain defined by an approximate solution. The content  of the paper was also presented in the book \cite{Agarwal2} two years later. It should be noticed that afer the appearance of this book many researchers studied  analytical solution \cite{Erturk}, \cite{Noor} and numerical solution \cite{Doedel}, \cite{Mohanty} of the problem \eqref{eq0} assuming the existence and uniqueness of a solution or referred to the book for the qualitative aspects of the problem.\par  

Motivated by the above facts, in this paper we use a novel method for existence and uniqueness and an iterative method for the problem \eqref{eq0}
which can overcome some limitations of \cite{Agarwal}. For simplicity of presentation we consider the problem 
\begin{equation}\label{eq1}
\begin{aligned}
u^{(4)}(x)&=f(x,u(x),u'(x),u''(x),u'''(x)), \quad 0 < x < 1, \\
u(0)&=u(1)=0,\quad u'(0)=u'(1)=0,
\end{aligned}
\end{equation}
which describes the deformations of an elastic beam with both fixed end-point. 

\par Differently from the approaches of the authors mentioned above for the problems for the fully fourth order equation, in this paper we use the method developed by ourselves recently in the works \cite{ALQ,AQ,Dang3,Dang4}. Namely, we reduce the problem \eqref{eq1} to an operator equation for a triplet of an unknown  function and two unknown numbers. With the assumption of continuity and boundnedness of the function $f(x,u,y,v,z)$ in a bounded domain we prove the existence of a solution. The uniqueness of the solution is established  under the additional assumption of Lipschitz condition. An iterative method for finding the solution is investigated. Many examples show the applicability of the obtained theoretical results and the efficiency of the iterative method. Especially, some of the examples show the advantage of our results over ones of \cite{Agarwal}.

\section{The existence and uniqueness of a solution}
To investigate the problem \eqref{eq1}, for $u\in C^4[0, 1]$ we set 
\begin{equation}\label{eq1.1}
\varphi (x)=f(x,u(x),u'(x),u''(x),u'''(x)). 
\end{equation}
Then problem becomes
\begin{equation*}\label{eq1.2*}
\begin{aligned}
&u^{(4)}(x)=\varphi (x), \quad 0 < x < 1, \\ 
&u(0)=u(1)=0,\quad u'(0)=u'(1)=0.
\end{aligned}
\end{equation*}
It has a unique solution
\begin{equation}\label{eq1.3*}
u(x)=\int_0^1 G_0(x,t)\varphi (t)dt,
\end{equation}
where $G_0(x,t)$ is the Green function 
\begin{equation*}\label{eq1.4*}
\begin{aligned}
G_0(x,t)=\dfrac{1}{6}\left\{\begin{array}{ll}
t^2(x-1)^2(3x-2tx-t), \quad 0\le t \le x \le 1,\\
x^2(t-1)^2(3t-2tx-x), \quad 0 \le x < t \le 1 .\\
\end{array}\right.
\end{aligned}
\end{equation*}
From \eqref{eq1.3*} we have
\begin{equation}\label{eq1.5*}
u'(x)=\int_0^1 {G_0'}_x(x,t)\varphi (t)dt,
\end{equation}
where
\begin{equation*}\label{eq1.6*}
{G_0'}_x(x,t)=\dfrac{1}{2}\left\{\begin{array}{ll}
t^2(1-x)(2tx-3x+1),  \quad 0\le t \le x \le 1,\\
(t-1)^2x(2t-2tx-x),  \quad 0 \le x < t \le 1 
\end{array}\right. 
\end{equation*}
It is easy verify that
\begin{equation}\label{eq1.7*}
\begin{aligned}
\int_0^1 |G_0(x,t)|dt \leq  \dfrac{1}{384}, \quad 
\int_0^1 |{G_0'}_x(x,t)|dt \leq  \dfrac{1}{72 \sqrt 3}, \quad \forall x \in [0,1].
\end{aligned}
\end{equation}
Now we set 
\begin{equation}\label{eq1.2}
\quad v(x)=u''(x), \quad v (0)=\alpha , \quad v(1)=\beta.
\end{equation}
Then the problem \eqref{eq1} is reduced to the following second order problems
\begin{equation}\label{eq5}
\left\{
\begin{array}{ll}
v''(x) =\varphi (x), \quad 0 < x < 1, \\ 
v(0)=\alpha , \quad v(1)=\beta ,
\end{array}
\right. 
\end{equation}
\begin{equation}\label{eq6}
\left\{
\begin{array}{ll}
u''(x)=v(x), \quad 0 < x < 1, \\ 
u(0)=u(1)=0.
\end{array}
\right. 
\end{equation}
It is obvious that $v=v_{\varphi,\alpha ,\beta } (x), u = u_{\varphi,\alpha ,\beta } (x)$, so from \eqref{eq1.1} and the conditions $u'(0)=u'(1)=0$ we have
\begin{equation}\label{eq7}
\left\{
\begin{array}{ll}
\varphi (x)=f(x,u_{\varphi,\alpha ,\beta } (x),u'_{\varphi,\alpha ,\beta } (x),u''_{\varphi,\alpha ,\beta } (x), u'''_{\varphi,\alpha ,\beta } (x)), \\ 
u'_{\varphi,\alpha ,\beta } (0)=0, \quad u'_{\varphi,\alpha ,\beta } (1)=0.
\end{array}
\right. 
\end{equation}
We notice that the solutions of the problem \eqref{eq5}, \eqref{eq6} can be represented in the forms
\begin{equation} \label{eqv}
v(x)=-\int_0^1 G(x,t)\varphi (t)dt+(\beta -\alpha )x+\alpha,
\end{equation}
\begin{equation*} \label{equ}
\begin{aligned}
u(x)&=-\int_0^1 G(x,s)v (s)ds\\
&=\int_0^1 G(x,s)\Big( \int_0^1 G(s,t)\varphi (t)dt\Big)ds-\int_0^1 G(x,s)\Big[(\beta -\alpha )s+\alpha \Big]ds,
\end{aligned}
\end{equation*}
where $G$ is the Green's function  
\begin{equation}\label{eqG}
\begin{aligned}
G(x,t)=\left\{\begin{array}{ll}
t(1-x),& 0\le t \le x \le 1,\\
x(1-t),&0 \le x \le t \le 1. \\
\end{array}\right.
\end{aligned}
\end{equation}
Therefore
\begin{equation*} \label{equ'}
\begin{aligned}
u'(x)&=\int_0^1 G'_x(x,s)\Big( \int_0^1 G(s,t)\varphi (t)dt\Big)ds-\int_0^1 G'_x(x,s)\Big[(\beta -\alpha )s+\alpha \Big]ds\\
&=-\int_0^x s\Big( \int_0^1 G(s,t)\varphi (t)dt\Big)ds+\int_x^1 (1-s)\Big( \int_0^1 G(s,t)\varphi (t)dt\Big)ds\\
&+\int_0^x s\Big[(\beta -\alpha )s+\alpha \Big]ds-\int_x^1 (1-s)\Big[(\beta -\alpha )s+\alpha \Big]ds.
\end{aligned}
\end{equation*}
By setting
\begin{equation} \label{eqH}
\begin{aligned}
H_0(t)&=\int_0^1 (1-s)G(s,t)ds=\dfrac{t^3}{6}-\dfrac{t^2}{2}+\dfrac{t}{3}, \\
H_1(t)&=\int_0^1 sG(s,t)ds=-\dfrac{t^3}{6}+\dfrac{t}{6}
\end{aligned}
\end{equation}
we obtain
\begin{equation} \label{equ'}
u'(0)=\int_0^1 H_0(t) \varphi (t)dt-\Big(\dfrac{\beta }{6}+\dfrac{\alpha }{3} \Big), \quad u'(1)=-\int_0^1 H_1(t) \varphi (t)dt +\Big(\dfrac{\beta }{3}+\dfrac{\alpha }{6} \Big).
\end{equation}
Combining \eqref{eq7} and \eqref{equ'} we have the system for unknowns $\varphi , \alpha , \beta $
\begin{equation}\label{eq8}
\left\{
\begin{array}{ll}
\varphi (x)=f(x,u_{\varphi,\alpha ,\beta } (x),u'_{\varphi,\alpha ,\beta } (x),u''_{\varphi,\alpha ,\beta } (x), u'''_{\varphi,\alpha ,\beta } (x)), \\ 
\displaystyle \int_0^1 H_0(t) \varphi (t)dt-\Big(\dfrac{\beta }{6}+\dfrac{\alpha }{3} \Big)=0,\\
-\displaystyle \int_0^1 H_1(t) \varphi (t)dt +\Big(\dfrac{\beta }{3}+\dfrac{\alpha }{6} \Big)=0.
\end{array}
\right. 
\end{equation}
\par Now we consider the Banach space $S=C[0,1]\times \mathbb R \times \mathbb R$ with the norm defined by
\begin{align}\label{norm}
\|\omega  \| =\| g \| _\infty +|a |+|b |, 
\end{align}
where
$\omega =(g \quad a \quad b  )^T \in S,$ $\| g \|_\infty =\max_{x \in [0,1]} | g(x)|.$
\par From \eqref{eq8}, it is easy to see that with $\omega =(\varphi  \quad \alpha  \quad \beta )^T \in S$, $$\omega =A\omega $$
where $A: S\rightarrow S$ be a nonlinear operator defined by
\begin{align}\label{eqA}
A\omega =\begin{pmatrix}
f(x,u_{\varphi,\alpha ,\beta } (x),y_{\varphi,\alpha ,\beta } (x),v_{\varphi,\alpha ,\beta } (x), z_{\varphi,\alpha ,\beta } (x)\\
3\displaystyle \int_0^1 H_0(t) \varphi (t)dt-\dfrac{\beta }{2}\\
3\displaystyle \int_0^1 H_1(t) \varphi (t)dt-\dfrac{\alpha  }{2}
\end{pmatrix},
\end{align}
with
\begin{equation}\label{eq9}
y_{\varphi,\alpha ,\beta } (x)=u'_{\varphi,\alpha ,\beta } (x), \quad z_{\varphi,\alpha ,\beta } (x)=v'_{\varphi,\alpha ,\beta } (x).
\end{equation} 
We will prove that under some conditions $A$ is a contraction operator.
For any $M>0,$ we define the set
\begin{equation} \label {eqD}
\begin{split}
\mathcal{D}_M= \Big \{ (x,u,y,v,z) \mid  0\leq x\leq 1; |u|\leq \dfrac{M}{384}; |y|\leq \dfrac{M}{72\sqrt 3};
|v| \leq M;|z|\leq M \Big \}.
\end{split}
\end{equation}
As usual, denote by $B[O,M]$ the closed ball centered at $O$ with radius $M$ in $S$. Notice that
\begin{align}\label{eqGG'}
\int_0^1 |G(x,t)|dt \leq \dfrac{1}{8}, \quad \forall x \in [0,1],
\end{align}
\begin{align}\label{eqH0H1}
\int_0^1 |H_0(t)|dt \leq \dfrac{21}{500}, \quad \int_0^1 |H_1(t)|dt \leq \dfrac{21}{500},
\end{align}
where the Green function $G(x,t)$ and the functions $H_0(t),$ $H_1(t)$ are defined by \eqref{eqG}, \eqref{eqH}, respectively.
\begin{Lemma}\label{lem1}
Suppose that the function $f$ is continuous  and there exists constant $M>0$ such that 
\begin{equation}\label{eq10'}
|f(x,u,y,v,z)| \leq \dfrac{M}{2}
\end{equation}
\noindent for all $(x,u,y,v,z) \in \mathcal{D}_M$. Then, the problem \eqref{eq1} has at least a solution.
\end{Lemma}
\noindent {\bf Proof.}
First we show that $A$ maps the closed ball $B[O,M]$ into itself. Indeed, let $\omega =(\varphi \quad \alpha \quad \beta )^T$ be an element  in $B[O,M]$. Then from the definition of the norm in \eqref{norm} we have
\begin{align}\label{eq10}
\| \omega \| =\| \varphi \| _\infty +|\alpha |+|\beta |\leq M.
\end{align}
We shall estimate $\| A\omega\|$ as follows: From \eqref{eq1.3*}-\eqref{eq1.7*} and \eqref{eq10} we get
\begin{align}\label{u}
\| u\| _\infty \leq \dfrac{1}{384}\| \varphi \|_\infty  \leq \dfrac{M}{384},
\end{align} 
\begin{align}\label{u'}
\| u'\| _\infty \leq \dfrac{1}{72\sqrt 3}\| \varphi \|_\infty  \leq \dfrac{M}{72\sqrt 3}.
\end{align} 
From \eqref{eq1.2}, \eqref{eqv}, \eqref{eqGG'} and \eqref{eq10} we obtain
\begin{align}\label{v}
\| u''\| _\infty=\| v\|_\infty  \leq \dfrac{1}{8}\| \varphi \| _\infty +\max(|\alpha |,|\beta |)\leq M,
\end{align} 
On the other hand, from \eqref{eqv} we have
\begin{equation}\label{equ'''}
\begin{aligned}
u'''(x)=v'(x)&=\Big[ \int_0^x t(x-1)\varphi (t)dt+\int_x^1 x(t-1)\varphi (t)dt +(\beta -\alpha )x\Big]'_x\\
&=\int_0^x t\varphi (t)dt+ \int_x^1 (t-1)\varphi (t)dt + \beta -\alpha .
\end{aligned}
\end{equation}
Therefore
 \begin{align}\label{u'''}
 \| u'''\| _\infty \leq \dfrac{1}{2}\| \varphi \| _\infty +|\beta-\alpha |\leq M,
 \end{align} 
Therefore, combining \eqref{u}-\eqref{v}, \eqref{u'''} with \eqref{eq9} and \eqref{eqD} we have 
\begin{align}\label{eq11}
(x,u,y,v,z) \in \mathcal{D}_M.
\end{align} 
On the other hand, from \eqref{eqH0H1} and \eqref{eq10} we obtain
\begin{equation}\label{eq12}
\begin{aligned}
&\Big | 3\displaystyle \int_0^1 H_0(t) \varphi (t)dt-\dfrac{\beta }{2}\Big | 
+ \Big | 3\displaystyle \int_0^1 H_1(t) \varphi (t)dt-\dfrac{\alpha  }{2}\Big |\\
&\leq \dfrac{63}{500}\| \varphi \| _\infty + \dfrac{\beta }{2}
+ \dfrac{63}{500}\| \varphi \| _\infty + \dfrac{\alpha  }{2}\leq \dfrac{M}{2}.
\end{aligned} 
\end{equation}
Hence, taking into account \eqref{eq10'}, \eqref{eq11} and \eqref{eq12} from the definition of operator $A$ in \eqref{eqA} we have
\begin{align*}
\| A\omega \| \leq M,
\end{align*} 
i.e., the operator $A$ maps $B[O,M]$ into itself.
\par Second, we shall show that $A$ is completely continuous in the ball $B[O,M]\subset S$. Indeed, suppose that for any $\omega =(\varphi \quad \alpha  \quad \beta )^T \in B[O,M]$, denote $u_\omega =u_{\varphi ,\alpha ,\beta }$ be the solution of the problem \eqref{eq1} for $\omega $ and also denote $u'_\omega=u'_{\varphi ,\alpha ,\beta }$, $u''_\omega=u''_{\varphi ,\alpha ,\beta }$, $u'''_\omega=u'''_{\varphi ,\alpha ,\beta }$. From the presentations of $u, u', u'', u'''$ in \eqref{eq1.3*}, \eqref{eq1.5*}, \eqref{eqv} and \eqref{equ'''} we infer that $u , u', u'', u'''$ are completely continuous from $B[O,M]$ to $C[0,1]$. 
\par Next, we will prove that for any bounded set $B_1 \subset B[O,M]$, the set $AB_1$ is relatively compact in $S$. Indeed, for any sequence $\{\omega _n\}=\{(\varphi_n \quad \alpha_n  \quad \beta_n )^T\} \in B_1$, from the complete continuity of 
$u$, there exist subsequence  
$$\{\omega^{(1)} _n\}=\{(\varphi^{(1)}_n \quad \alpha^{(1)}_n  \quad \beta^{(1)}_n )^T\}\subset \{\omega _n\}$$ such that
$$\lim_{n\rightarrow \infty }u_{\omega^{(1)} _n}=u^{(0)} \in C[0,1].$$ 
From the complete continuity of $u'$, there exist subsequence  
$$\{\omega^{(2)} _n\}=\{(\varphi^{(2)}_n \quad \alpha^{(2)}_n  \quad \beta^{(2)}_n )^T\}\subset \{\omega^{(1)} _n\}$$ such that
$$\lim_{n\rightarrow \infty }u'_{\omega^{(2)} _n}=y^{(0)} \in C[0,1].$$ 
From the complete continuity of $u''$, there exist subsequence  
$$\{\omega^{(3)} _n\}=\{(\varphi^{(3)}_n \quad \alpha^{(3)}_n  \quad \beta^{(3)}_n )^T\}\subset \{\omega^{(2)} _n\}$$ such that
$$\lim_{n\rightarrow \infty }u''_{\omega^{(3)} _n}=v^{(0)} \in C[0,1].$$
From the complete continuity of $u'''$, there exist subsequence  
$$\{\omega^{(4)} _n\}=\{(\varphi^{(4)}_n \quad \alpha^{(4)}_n  \quad \beta^{(4)}_n )^T\}\subset \{\omega^{(3)} _n\}$$ such that
$$\lim_{n\rightarrow \infty }u'''_{\omega^{(4)} _n}=z^{(0)} \in C[0,1].$$ 
Therefore, from the continuity of function $f$ we have 
$$\lim_{n\rightarrow \infty }f(x,u_{\omega^{(4)} _n},y_{\omega^{(4)} _n},v_{\omega^{(4)} _n},z_{\omega^{(4)} _n})
=f(x,u^{(0)},y^{(0)},v^{(0)},z^{(0)})\in C[0,1].$$ 
Further, it is possible to extract a subsequence  
$$\{\omega^{(5)} _n\}=\{(\varphi^{(5)}_n \quad \alpha^{(5)}_n  \quad \beta^{(5)}_n )^T\}\subset \{\omega^{(4)} _n\}$$ such that
$$\lim_{n\rightarrow \infty }\Big ( 3\displaystyle \int_0^1 H_0(t) \varphi^{(5)}_n(t)dt-\dfrac{\beta^{(5)}_n }{2} \Big )=a^{(0)} \in \mathbb R$$ 
and a subsequence  
$$\{\omega^{(6)} _n\}=\{(\varphi^{(6)}_n \quad \alpha^{(6)}_n  \quad \beta^{(6)}_n )^T\}\subset \{\omega^{(5)} _n\}$$ such that
$$\lim_{n\rightarrow \infty }\Big ( 3\displaystyle \int_0^1 H_1(t) \varphi^{(6)}_n(t)dt-\dfrac{\alpha^{(6)}_n}{2}\Big)=b^{(0)} \in \mathbb R.$$ 
Therefore, for any sequence $\{\omega _n\}=\{(\varphi_n \quad \alpha_n  \quad \beta_n )^T\} \in B_1$, there exist subsequence 
$$\{\omega^{(6)} _n\}=\{(\varphi^{(6)}_n \quad \alpha^{(6)}_n  \quad \beta^{(6)}_n )^T\}\subset \{\omega_n\}$$ such that
\begin{align*}
\lim_{n\rightarrow \infty }A \omega^{(6)} _n =\begin{pmatrix}
f(x,u^{(0)},y^{(0)},v^{(0)}, z^{(0)})\\
a^{(0)}\\
b^{(0)}
\end{pmatrix} \in S.
\end{align*}
Hence, $A$ is completely continuous in the ball $B[O,M]\subset S$. By the Schauder's fixed point theorem, $A$ has at least a fixed point, i.e., the problem \eqref{eq1} has at least one solution. The lemma is proved. 
\begin{Theorem}\label{theo1}
Suppose that the assumptions of Lemma \ref{lem1} hold. Further assume that there exists constants $K_1, K_2, K_3, K_4 \geq 0$ such that
\begin{equation}\label{eq11'}
\begin{split}
|f(x,u_2,y_2,v_2,z_2)-f(x,u_1,y_1,v_1,z_1)| & \leq K_1 |u_2-u_1|+K_2|y_2-y_1|\\
&+K_3|v_2-v_1|+K_4|z_2-z_1|,
\end{split}
\end{equation}
\noindent for all $(x,u_i,y_i,v_i,z_i) \in \mathcal{D}_M$ $(i=1,2)$ and
\begin{equation} \label {eq12'}
\begin{split}
q= \dfrac{K_1}{384} + \dfrac{K_2}{72\sqrt 3} + K_3 + K_4 < \dfrac{1}{2}
\end{split}
\end{equation}
then the problem \eqref{eq1} has a unique solution $u$ and
\begin{equation*}\label{eq1.16}
\| u\|_\infty \leq \dfrac{M}{384}, \quad \| u'\|_\infty \leq \dfrac{M}{72\sqrt 3}, 
\quad \| u''\|_\infty \leq M, \quad \| u'''\|_\infty \leq M.
\end{equation*}
\end{Theorem}
\noindent {\bf Proof.}
\par For proving the theorem, we shall show that the operator $A$ mentioned above is a contraction operator tin the ball $B[O,M]$. Indeed, suppose that $\omega  _1 =(\varphi _1 \quad \alpha _1 \quad \beta _1)^T$, $\omega  _2 =(\varphi _2 \quad \alpha _2 \quad \beta _2)^T \in B[O,M]$. Denote by $u_1,u_2$  the solutions of the problems \eqref{eq5}, \eqref{eq6}, respectively. We also denote $y_i=u_i',$ $v_i=u_i'',$ $z_i=u_i'''$ $(i=1,2).$ Then, as induced above $(x,u_i,y_i,v_i,z_i) \in \mathcal{D}_M$ $(i=1,2).$ Due to the estimates \eqref{u} -\eqref{u'''} and \eqref{eq12} we obtain
\begin{equation}\label{eq13}
\begin{aligned}
\| u_2-u_1\|_\infty   &\leq \dfrac{1}{384}\| \omega _2-\omega _1 \|, \quad \| y_2-y_1\|_\infty  \leq \dfrac{1}{72\sqrt 3}\| \omega _2-\omega  _1 \|, \\
\| v_2-v_1\|_\infty  &\leq \| \omega _2-\omega _1 \|, \quad \| z_2-z_1\|_\infty  \leq \| \omega _2-\omega _1 \|,
\end{aligned}
\end{equation}
\begin{equation}\label{eq14}
\begin{aligned}
&\Big | 3\displaystyle \int_0^1 H_0(t) (\varphi_2 (t)-\varphi_1 (t))dt-\dfrac{\beta_2-\beta_1 }{2}\Big | \\
&+ \Big | 3\displaystyle \int_0^1 H_1(t) (\varphi_2 (t)-\varphi_1 (t))dt-\dfrac{\alpha_2-\alpha _1 }{2}\Big |\leq \dfrac{\| \omega _2 -\omega _1\| }{2}.
\end{aligned}
\end{equation}
From \eqref{eq11'} and \eqref{eq13} we have
\begin{equation}\label{eq15}
\begin{aligned}
&|f(x,u_2,y_2,v_2,z_2)-f(x,u_1,y_1,v_1,z_1)|\\
&\leq K_1 |u_2 -u_1|+K_2 |y_2 -y_1|+K_3 |v_2 -v_1|+K_4 |z_2 -z_1|\\
&\leq \Big( \dfrac{K_1}{384}+\dfrac{K_2}{72\sqrt 3}+K_3+K_4 \Big)\| \omega _2 - \omega _1\| .
\end{aligned}
\end{equation}
Combining \eqref{eq14} and \eqref{eq15} we obtain
\begin{equation*}
\|A\omega _2-A\omega _1\| \leq \Big( \dfrac{K_1}{384}+\dfrac{K_2}{72\sqrt 3}+K_3+K_4 +\dfrac{1}{2}\Big)\| \omega _2-\omega _1\|, 
\end{equation*}
Therefore, taking into account \eqref {eq12'} we conclude that $A$ is a contraction operator in $B[O,M].$ Therefore, the equation $\omega  = A \omega  $ has a unique solution $\omega $  with $\|\omega  \| \le M$. This fact implies that the problem \eqref{eq1} has a unique solution $u$ determined from \eqref{eq5}-\eqref{eq6} with the found triplet $\omega =(\varphi \quad \alpha \quad \beta )^T$. The estimates for $u, u', u'', u'''$ follow straightforward from the estimates \eqref{u}-\eqref{u'''}. The theorem is proved.
\section{Solution method and numerical examples}
Consider the following iterative method for solving the problem \eqref{eq1}:\\
\noindent 
i) Given an initial approximation $\varphi _0(x), \alpha _0, \beta _0$, for example, 
\begin{equation}\label{eq32}
\varphi _0(x)=f(x,0,0,0,0), \quad \alpha _0=0, \quad \beta _0=0.
\end{equation}
ii) Knowing $\varphi _k(x), \alpha _k, \beta _k$ $(k=0,1,2,...)$ solve consecutively two problems
\begin{equation}\label{eq33}
\left\{
\begin{array}{ll}
v_k''(x) =\varphi_k (x), \quad 0 < x < 1, \\ 
v_k(0)=\alpha _k, \quad v_k(1)=\beta _k,
\end{array}
\right.
\end{equation}
\begin{equation}\label{eq34}
\left\{
\begin{array}{ll}
u_k''(x) =v_k (x), \quad 0 < x < 1, \\ 
u_k(0)=u_k(1)=0.
\end{array}
\right.\end{equation}
iv) Update the new approximation
\begin{equation}\label{eq35}
\begin{aligned}
\varphi_{k+1}(x)&=f(x,u_k(x),u'_k(x),u''_k(x),u'''_k(x)),\\
\alpha_{k+1}&=3\displaystyle \int_0^1 H_0(t) \varphi_{k+1} (t)dt-\dfrac{\beta _k  }{2},\\
\beta _{k+1}&=3\displaystyle \int_0^1 H_1(t) \varphi_{k+1} (t)dt-\dfrac{\alpha_{k+1}}{2}.
\end{aligned}
\end{equation}
Set $p_k=\dfrac{(q+\frac{1}{2})^k}{\frac{1}{2}-q}\| \omega _1-\omega _0\| .$ We have the following result
\begin{Theorem}\label{theo2}
Under the assumptions of Theorem \ref{theo1}, the above iterative method converges with the rate of geometric progression and there hold the estimates 
\begin{equation*} \label{eq1.23}
\begin{split}
\| u_k-u\|_\infty &\leq \dfrac{p_k}{384}, \quad \| u'_k-u'\|_\infty \leq \dfrac{p_k}{72\sqrt 3},\\
\| u''_k-u''\|_\infty &\leq p_k, \quad \| u'''_k-u'''\|_\infty \leq p_k,
\end{split}
\end{equation*}
 where $u$ is the exact solution of the problem \eqref{eq1}.
\end{Theorem} 
\noindent {\bf Proof.}
Notice that the above iterative method is the successive iteration method for finding the fixed point of the operator
$A$ with the initial approximation  $\omega_0=(\varphi _0(x) \quad \alpha _0 \quad \beta _0)^T \in B[O,M]$. Therefore, it converges with the rate of geometric progression and there is the estimate
\begin{equation*}
||\omega  _k - \omega  || \leq \frac{(q+\frac{1}{2})^k}{\frac{1}{2}-q}||\omega  _1 - \omega  _0 ||.
\end{equation*}
Combining it with the estimate of the type \eqref{eq13} we obtain the results of the theorem.\\

Below we illustrate the obtained theoretical results of the existence and uniqueness of a solution and the convergence of the iterative method on some examples, where the exact solution of the problem is known or is not known.
\par In order to numerically realize the iterative process we use the difference scheme of fourth order accuracy  \cite{Sam} for solving \eqref{eq33}, \eqref{eq34}, and formulas of the same order of accuracy for approximating
the first derivative, the second derivative, the third derivative and the definite integral on the uniform grid
$ \overline \omega_h=\{x_i=ih,\quad i=0,1,...,N\} .$  For grid functions we use the uniform norm defined as $\| u\| _{\overline \omega_h}=\max_{x_i \in \overline \omega_h} |u(x_i)|.$ \par

For  testing the convergence of the proposed iterative method we perform some experiments for the case of known exact solutions and also for the case of unknown exact solutions. 

In all the tables and the graphs of computation below, $u$ is the exact solution, $N$ denotes the number of grid intervals, $K$ denotes the number of performed iterations, $eu(k)=\| u_k-u\|_{_{\overline \omega_h}}$, $e(k)=\| u_k-u_{k-1}\|_{_{\overline \omega_h}}$. We perform the iterative process until $e(k) \le 10^{-15} $. \\
\par
\begin{figure}
\begin{center}
\includegraphics[height=6cm,width=11cm]{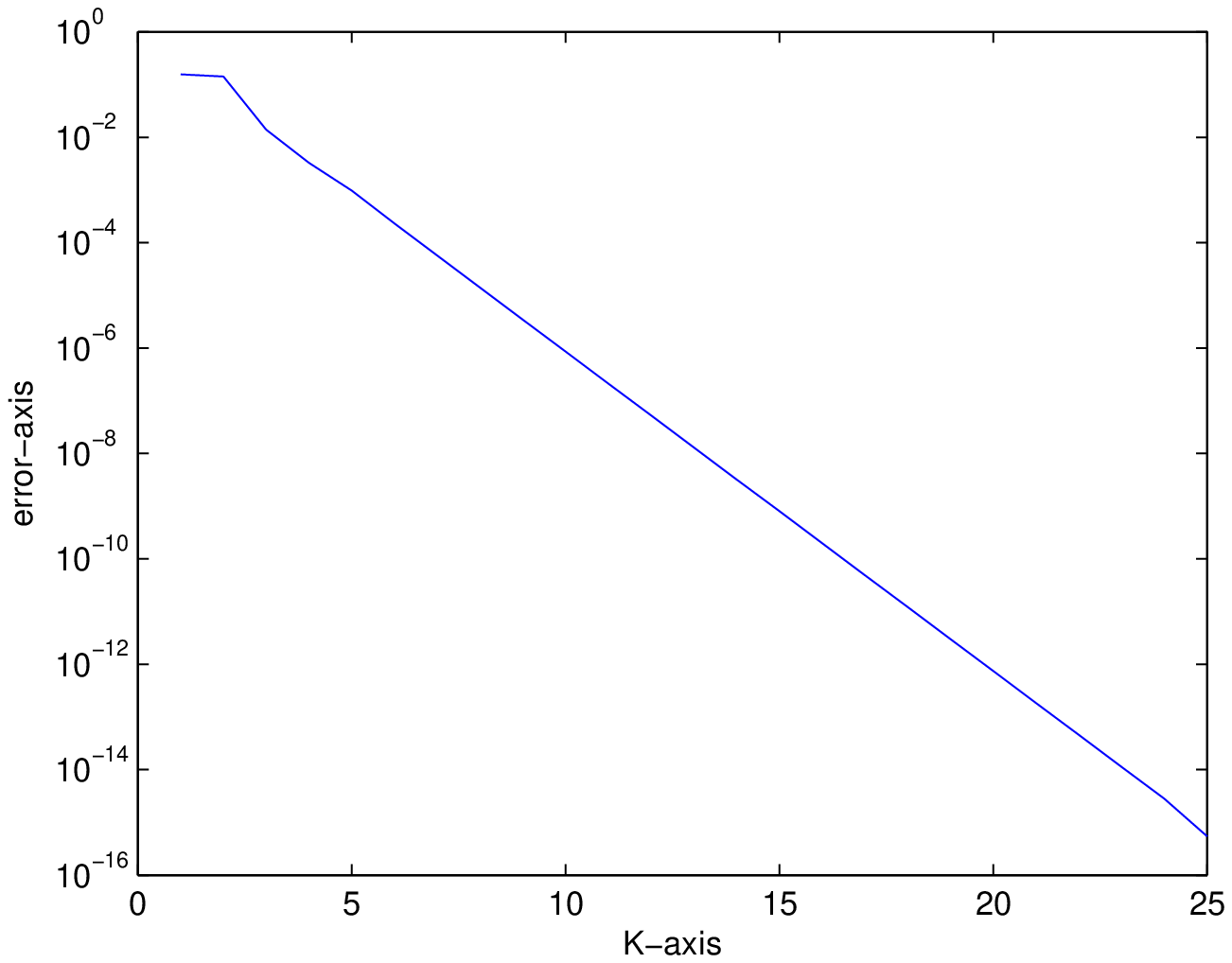}
\caption{The graph of $e(k)$ in Example $1$. }
\label{fig1.1}
\end{center}
\end{figure}

 First, we consider the example for the case of known exact solution.\\
{\bf{Example 1.} } Consider the problem
 \begin{equation*}
 \left\{
\begin{array}{lll}
u^{(4)}(x)=12+\dfrac{u(x)u'''(x)}{2}-\dfrac{u'(x)u''(x)}{4}+\dfrac{u'}{4}, \quad 0 < x <1,\\
u(0)=u(1)=0, \quad u'(0)=u'(1)=0.
\end{array}\right.
\end{equation*}
The exact solution of the problem is
 $$u(x)=\dfrac{x^4}{2}-x^3+\dfrac{x^2}{2}, \quad 0\leq x\leq 1.$$
We have 
$$f(x,u,y,v,z)=12+\dfrac{uz}{2}-\dfrac{yv}{4}+\dfrac{y}{4}.$$
In the domain 
\begin{equation*} 
\begin{split}
\mathcal{D}_M= \Big \{ (x,u,y,v,z) \mid  0\leq x\leq 1; |u|&\leq \dfrac{M}{384} ; |y|\leq \dfrac{M}{72\sqrt 3} ;
|v| \leq M ;|z|\leq M  \Big \}
\end{split}
\end{equation*}
we have
\begin{equation*}
 \begin{aligned}
 |f(x,u,y,v,z)|&\leq 12+\dfrac{1}{2}\dfrac{M}{384}.M+\dfrac{1}{4}.\dfrac{M}{72\sqrt 3}.M+\dfrac{1}{4}.\dfrac{M}{72\sqrt 3}\\
 &=12+\dfrac{M^2}{768}+\dfrac{M^2}{288\sqrt 3}+\dfrac{M}{288\sqrt 3}.
 \end{aligned}
 \end{equation*}
Therefore, the choice $M=36$ guarantees the satisfaction of the condition 
\eqref{eq10'}. Besides, in the domain $\mathcal{D}_{36}$, since
\begin{equation*}
\begin{aligned}
|f'_u|=\dfrac{|z|}{2}&\leq \dfrac{M}{2}, \quad |f'_y| = \dfrac{|v|}{4}+\dfrac{1}{4}\leq \dfrac{M}{4}+\dfrac{1}{4}, \\
\quad |f'_v| &= \dfrac{|y|}{4}\leq \dfrac{M}{288\sqrt 3}, 
\quad |f'_z| = \dfrac{|u|}{2}\leq \dfrac{M}{768},
\end{aligned}
\end{equation*} 
we can take 
 $$K_1=18, \quad K_2=\dfrac{37}{4}, \quad K_3=\dfrac{1}{8\sqrt 3}, \quad K_4=\dfrac{3}{64}.$$
Then
 $$q= \dfrac{K_1}{384} + \dfrac{K_2}{72\sqrt 3} + K_3+ K_4 \approx  0.24<\dfrac{1}{2}.$$
All the conditions in Theorem \ref{theo1} are satisfied. Consequently, the problem has a unique solution, and by Theorem \ref{theo2} the iterative method \eqref{eq32}-\eqref{eq35} converges.
Table 1 shows the convergence of the iterative method. From the table we see that the convergence of the discrete version of the iterative method does not depend on the grid size.\\
\begin{table}[!ht]
\centering
\setlength{\tabcolsep}{0.6cm}
\caption[smallcaption]{The convergence in Example 1 }
\label{Tab1}
\begin{tabular}{ cccc} 
\hline
$N$ & $K$ & $ eu(k)$ & $e(k)$ \\
\hline
100&25&2.8103e-11&8.0491e-16\\
200&25&3.7123e-12&6.8348e-16\\
500&25&6.5989e-13&6.9736e-16\\
1000&25&3.6911e-14&7.0777e-16\\
\hline
\end{tabular}
\end{table}
The graphs of $e(k)$ of Example 1 is depicted in Figure \ref{fig1.1}.\\
\par
\begin{figure}
\begin{center}
\includegraphics[height=6cm,width=11cm]{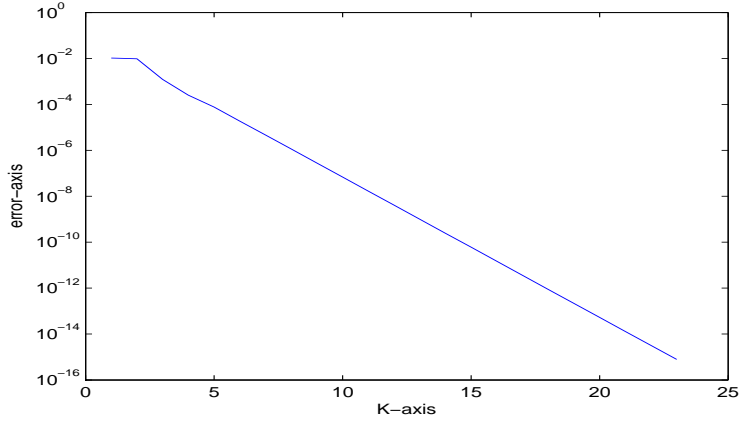}
\caption{The graph of $e(k)$ in Example $2$. }
\label{fig2.1}
\end{center}
\end{figure}

\begin{figure}
\begin{center}
\includegraphics[height=6cm,width=11cm]{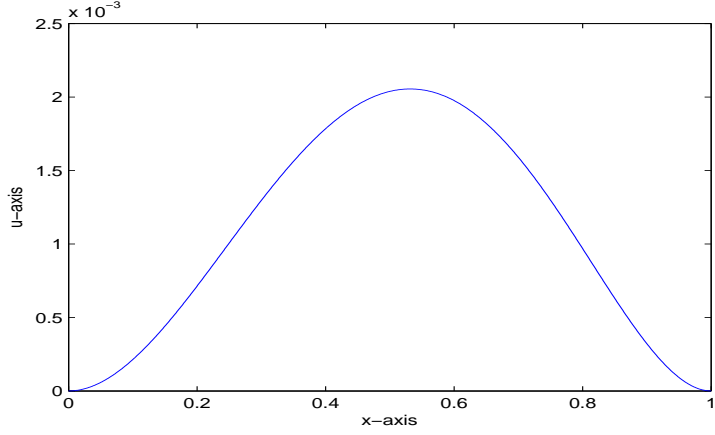}
\caption{The graph of the approximate solution in Example $2$. }
\label{fig2.2}
\end{center}
\end{figure}

 In the next example, the exact solution of the problem \eqref{eq1} is not known.\\
{\bf{Example 2.} } Consider the problem
 \begin{equation*}
 \left\{
\begin{array}{lll}
u^{(4)}(x)=x+x^2+u^2(x)u''(x)+u'(x) \sin u'''(x), \quad 0 < x <1,\\
u(0)=u(1)=0, \quad u'(0)=u'(1)=0.
\end{array}\right.
\end{equation*}
We have 
\begin{equation*} \label{Ex2.2}
f(x,u,y,v,z)=x+x^2+ u^2v+y\sin z.
\end{equation*}
In the domain 
\begin{equation*} 
\begin{split}
\mathcal{D}_M= \Big \{ (x,u,y,v,z) \mid  0\leq x\leq 1; |u|&\leq \dfrac{M}{384} ; |y|\leq \dfrac{M}{72\sqrt 3} ;
|v| \leq M ;|z|\leq M  \Big \}
\end{split}
\end{equation*}
we have
\begin{equation*}
|f(x,u,y,v,z)|\leq 2+\dfrac{M^3}{147456}+\dfrac{M}{72\sqrt 3}.
\end{equation*}
Therefore, the choice $M=5$ guarantees the satisfaction of the condition \eqref{eq10'}. Besides, in the domain $\mathcal{D}_{5}$, since
\begin{equation*}
\begin{aligned}
|f'_u|&=2|uv|\leq \dfrac{M^2}{192}, \quad |f'_y| = |\sin z|\leq 1, \\
\quad |f'_v| &= |u^2|\leq \dfrac{M^2}{147456}, 
\quad |f'_z| = |y\cos z|\leq \dfrac{M}{72\sqrt 3},
\end{aligned}
\end{equation*} 
we can take 
 $$K_1=\dfrac{25}{192}, \quad K_2=1, \quad K_3=\dfrac{25}{147456}, \quad K_4=\dfrac{5}{72\sqrt 3}.$$
Then
 $$q= \dfrac{K_1}{384} + \dfrac{K_2}{72\sqrt 3} + K_3+ K_4 \approx  0.05<\dfrac{1}{2}.$$
All the conditions in Theorem \ref{theo1} are satisfied. Consequently, the problem has a unique solution, and by Theorem \ref{theo2} the iterative method \eqref{eq32}-\eqref{eq35} converges. The results of computation show that for different grid sizes the discrete version of the iterative process reached the tolerance $e(k)=10^{-15}$ after 23 iterations. The graphs of $e(k)$ and the approximate solution are depicted in Figure \ref{fig2.1} and Figure \ref{fig2.2}, respectively.\\

\begin{figure}
\begin{center}
\includegraphics[height=6cm,width=11cm]{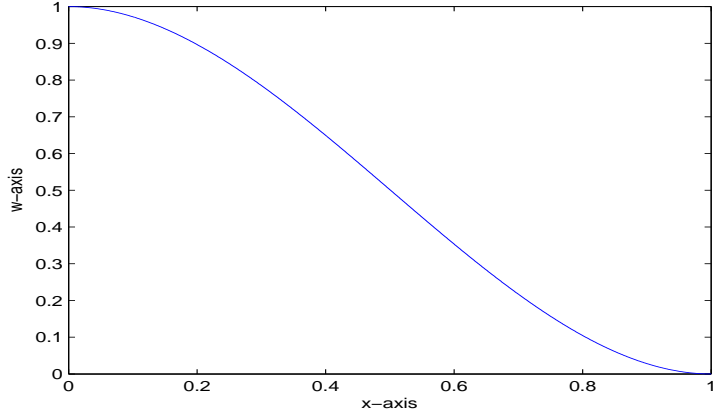}
\caption{The graph of the approximate solution in Example $3$. }
\label{fig3.1}
\end{center}
\end{figure}

{\bf{Example 3.} } Consider the example (see \cite[Example 7.1] {Agarwal})
 \begin{equation} \label{Exam3}
 \left\{
\begin{array}{lll}
w^{(4)}(x)=w^2(x)\sin w(x)+\sin x, \quad 0 < x <1,\\
w(0)=1, \quad w(1)=w'(0)=w'(1)=0.
\end{array}\right.
\end{equation}
We set $u(x)=w(x)-P(x),$ where $P(x)=2x^3-3x^2+1.$ Then the problem \eqref{Exam3} becomes 
\begin{equation} \label{Exam3'}
 \left\{
\begin{array}{lll}
u^{(4)}(x)=[u(x)+P(x)]^2\sin [u(x)+P(x)]+\sin x, \quad 0 < x <1,\\
u(0)=u(1)=0, \quad u'(0)=u'(1)=0.
\end{array}\right.
\end{equation}
We consider problem \eqref{Exam3'}. We have 
$$f(x,u,y,v,z)=(u+P)^2\sin [u+P]+\sin x.$$
In the domain 
\begin{equation*} 
\begin{split}
\mathcal{D}_M= \Big \{ (x,u,y,v,z) \mid  0\leq x\leq 1; |u|&\leq \dfrac{M}{384} ; |y|\leq \dfrac{M}{72\sqrt 3} ;
|v| \leq M ;|z|\leq M  \Big \}
\end{split}
\end{equation*}
we have
\begin{equation*}
|f(x,u,y,v,z)|\leq 2+\dfrac{M^2}{147456}+\dfrac{M}{192}.
\end{equation*}
Therefore, the choice $M=6$ guarantees the satisfaction of the condition \eqref{eq10'}. Besides, in the domain $\mathcal{D}_{6}$, since
\begin{equation*}
\begin{aligned}
|f'_u|=|2(u+P)\sin(u+P)&+(u+P)^2\cos (u+P)|\leq \dfrac{M^2}{147456}+\dfrac{M}{96}+3, \\
 |f'_y| &=|f'_v|=|f'_z| = 0,
\end{aligned}
\end{equation*} 
we can take 
 $$K_1=\dfrac{12545}{4096}, \quad K_2= K_3= K_4=0.$$
Then
 $$q= \dfrac{K_1}{384} + \dfrac{K_2}{72\sqrt 3} + K_3+ K_4 \approx  0.008<\dfrac{1}{2}.$$
All the conditions in Theorem \ref{theo1} are satisfied. Consequently, the problem \eqref{Exam3'} has a unique solution, and therefore, so does the problem \eqref{Exam3}.  Besides, by Theorem \ref{theo2} the iterative method \eqref{eq32}-\eqref{eq35} converges. The results of computation show that for different grid sizes the discrete version of the iterative process reached the tolerance $e(k)=10^{-15}$ after 23 iterations. The graphs of the approximate solution of the problem \eqref{Exam3} is depicted in Figure \ref{fig3.1}.
\par
\begin{figure}
\begin{center}
\includegraphics[height=6cm,width=11cm]{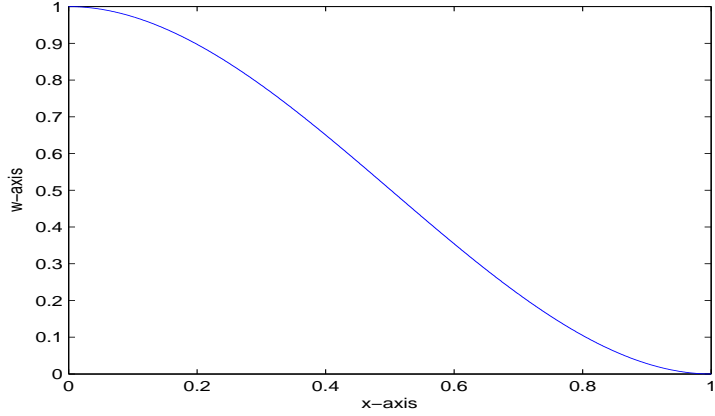}
\caption{The graph of the approximate solution in Example $4$. }
\label{fig4.1}
\end{center}
\end{figure}

 Remark that in \cite{Agarwal} the authors could only establish the existence of a solution of the problem \eqref{Exam3} in the region
$$S=\Big\{(x,w):0\leq x\leq 1,|w|\leq k_0 \quad \text{where} \quad 1\leq k_0\leq 383.9973958... \Big\}$$
but did not guarantee the uniqueness of a solution of this problem meanwhile by using Theorem \ref{theo1} in this paper, the problem \eqref{Exam3} has a unique solution in the region
$$S=\Big\{(x,w):0\leq x\leq 1,|w-P|\leq \dfrac{6}{384} \Big\}$$
i.e.,
$$S=\Big\{(x,w):0\leq x\leq 1,|w|\leq 1.015625 \Big\}.$$
{\bf{Example 4.} } Consider the example in \cite[Example 7.3] {Agarwal}
 \begin{equation} \label{Exam4}
 \left\{
\begin{array}{lll}
w^{(4)}(x)=w(x)\sin w(x)+e^{-x^2}, \quad 0 < x <1,\\
w(0)=1, \quad w(1)=w'(0)=w'(1)=0.
\end{array}\right.
\end{equation}
We set $u(x)=w(x)-P(x),$ where $P(x)=2x^3-3x^2+1.$ Then the problem \eqref{Exam4} becomes 
\begin{equation} \label{Exam4'}
 \left\{
\begin{array}{lll}
u^{(4)}(x)=[u(x)+P(x)]\sin [u(x)+P(x)]+e^{-x^2}, \quad 0 < x <1,\\
u(0)=u(1)=0, \quad u'(0)=u'(1)=0.
\end{array}\right.
\end{equation}
Analogously as in Example 3, we can choose $M=6$, and therefore, the Lipschitz coefficients in Theorem \ref{theo1} are $K_1=\dfrac{129}{64}, K_2=K_3=K_4=0$. Then $q \approx 0.005<\dfrac{1}{2}.$
All the conditions of Theorem \ref{theo1} are satisfied. Hence the problem \eqref{Exam4'} (therefore the problem \eqref{Exam4}) has a unique solution and the iterative method converges. The results of computation show that for different grid sizes the discrete version of the iterative process reached the tolerance $e(k)=10^{-15}$ after 24 iterations. The graph of the approximate solution of the problem \eqref{Exam4} is depicted in Figure \ref{fig4.1}. \\

\par Remark that in \cite{Agarwal} the authors can only establish the existence of a solution of the problem \eqref{Exam4} in the region
$$S=\Big\{(x,w):0\leq x\leq 1,|w|\leq 1.005222 \Big\}$$
but did not guarantee the uniqueness of a solution of this problem meanwhile by using Theorem \ref{theo1} in this paper, the problem \eqref{Exam4} has a unique solution in the region
$$S=\Big\{(x,w):0\leq x\leq 1,|w-P|\leq \dfrac{6}{384} \Big\}$$
i.e.,
$$S=\Big\{(x,w):0\leq x\leq 1,|w|\leq 1.015625 \Big\}.$$
{\bf{Example 5.} } Consider the  example in \cite[Example 7.2] {Agarwal}
 \begin{equation} \label{Exam5}
 \left\{
\begin{array}{lll}
w^{(4)}(t)=w(t)^{1/2}\sin e^{w(t)}+e^{-t^2}, \quad a < t <b,\\
w(a)=A_1, \quad w(b)=B_1, \quad w'(a)=A_2, \quad w'(b)=B_2.
\end{array}\right.
\end{equation}
We set $v(t)=w(t)-P(t),$ where $P(t)$ is the third degree polynomial satisfy the boundary conditions in the problem \eqref{Exam5}. Then this problem becomes 
\begin{equation} \label{Exam5'}
 \left\{
\begin{array}{lll}
v^{(4)}(t)=[v(t)+P(t)]^{1/2}\sin e^{v(t)+P(t)}+e^{-t^2}, \quad a < t <b,\\
v(a)=v(b)=0, \quad v'(a)=v'(b)=0.
\end{array}\right.
\end{equation}
Now, we set $t=a+(b-a)x$ and $u(x)=v(a+(b-a)x)$. Then the problem \eqref{Exam5'} becomes
\begin{equation} \label{Exam5''}
 \left\{
\begin{array}{lll}
u^{(4)}(x)=(b-a)^4\big\{ [u(x)+P(a+(b-a)x)]^{1/2}\sin e^{u(x)+P(a+(b-ax))}\\
+e^{-[a+(b-a)x]^2} \Big\}, \quad 0 < x <1,\\
u(0)=u(1)=0, \quad u'(0)=u'(1)=0.
\end{array}\right.
\end{equation}
We have 
$$f(x,u,y,v,z)=(b-a)^4\Big \{ [u+P(a+(b-a)x)]^{1/2}\sin e^{u+P(a+(b-a)x)}+e^{-[a+(b-a)x]^2} \Big \}.$$
In the domain 
\begin{equation*} 
\begin{split}
\mathcal{D}_M= \Big \{ (x,u,y,v,z) \mid  0\leq x\leq 1; |u|&\leq \dfrac{M}{384} ; |y|\leq \dfrac{M}{72\sqrt 3} ;
|v| \leq M ;|z|\leq M  \Big \}
\end{split}
\end{equation*}
we have
\begin{equation*}
|f(x,u,y,v,z)|\leq (b-a)^4 \Big\{ \sqrt{\dfrac{M}{384}+K} + e^{-\min (a^2, b^2)} \Big\},
\end{equation*}
where $K=\max_{0\leq x\leq 1}|P(a+(b-a)x)|.$
We can see that for any finite numbers $b-a, A_1, A_2, B_1, B_2$ there exists constant $M>0$ such that the condition \eqref{eq10'} is satisfied. Therefore, the problem \eqref{Exam5''} has at least one solution, i.e., the problem \eqref{Exam5} has at least one solution.
\par Remark that in \cite{Agarwal}, the authors also established the existence of a solution of the problem \eqref{Exam5}.
\par
\begin{figure}
\begin{center}
\includegraphics[height=6cm,width=11cm]{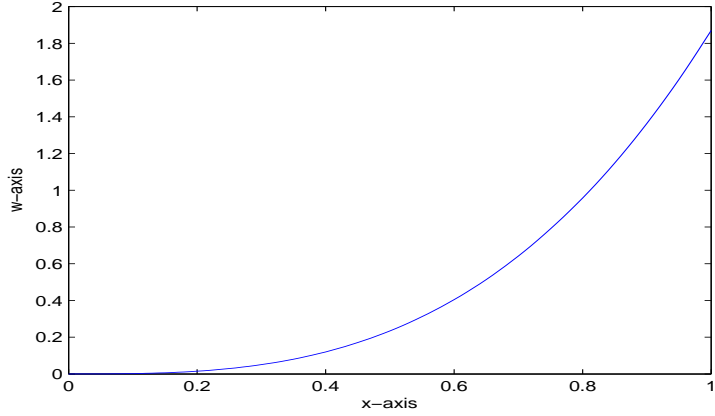}
\caption{The graph of the approximate solution in Example $6$. }
\label{fig6.1}
\end{center}
\end{figure}

 In the next example we shall show that the conditions in the theorems of the existence of solutions in \cite{Agarwal} are not satisfied meanwhile by using Theorem \ref{theo1} in this paper, the problem has a unique solution.\\
{\bf{Example 6.} } Consider the problem
  \begin{equation} \label{Exam6}
 \left\{
 \begin{array}{lll}
 w^{(4)}(x)=w^5(x), \quad 0 < x <1,\\
 w(0)=0, \quad w(1)=1.87, \quad w'(0)=0, \quad w'(1)=5.61.
 \end{array}\right.
 \end{equation}
It is easy to see that \cite[Corollary 3.3, Theorem 3.5]{Agarwal}) cannot be applied in this example. Next we shall show that the conditions in \cite[Theorem 3.1]{Agarwal}) are not satisfied. Indeed, we have $P_3(x)=1.87x^3$ (the third degree polynomial satisfying the boundary conditions in the problem \eqref{Exam6}). Let
$$\max_{0\leq x\leq 1}|P_3^{(i)}(x)|\leq k_i, i=0,1,2,3, k_i>0.$$
Clearly, $k_0\geq 1.87.$ Consider the compact set
$$D=\{ (x,w_0,w_1,w_2,w_3): 0\leq x\leq 1, |w_i|\leq 2k_i, i=0,1,2,3 \}.$$
It is obvious that
$$Q=\max_D|f(x,w_0,w_1,w_2,w_3)|=\max_D|w_0^5|=32k_0^5.$$
Therefore, since $C_{4,0}=\dfrac{1}{384}$ we have
$$\dfrac{k_0}{QC_{4,0}}=\dfrac{384k_0}{32k_0^5}=\dfrac{12}{k_0^4}<1=(1-0), \forall k_0\geq 1.87,$$
i.e., the condition (iii) in \cite[Theorem 3.1]{Agarwal}) is not satisfied.  
 Hence, the theorems of the existence of solutions in \cite{Agarwal} cannot be applied in this example. \\
Now, by using Theorem \ref{theo1}, we shall show that the problem has a unique solution. For this purpose we set 
$u(x)=w(x)-P(x),$ where $P(x)=1.87x^3.$ Then the problem \eqref{Exam6} becomes 
 \begin{equation} \label{Exam6'}
  \left\{
 \begin{array}{lll}
 u^{(4)}(x)=\Big(u(x)+P(x)\Big)^5, \quad 0 < x <1,\\
 u(0)=u(1)=0, \quad u'(0)=u'(1)=0.
 \end{array}\right.
 \end{equation}
 Analogously as in Example 3 and Example 4, we can choose $M=100$, and therefore, the Lipschitz coefficients in Theorem \ref{theo1} are $K_1=103, K_2=K_3=K_4=0$. Then $q \approx 0.27<\dfrac{1}{2}.$
 All the conditions of Theorem \ref{theo1} are satisfied. Hence the problem \eqref{Exam6'}, and in consequence, the problem \eqref{Exam6} has a unique solution and the iterative method converges. The results of computation show that for different grid sizes the discrete version of the iterative process reached the tolerance $e(k)=10^{-15}$ after 23 iterations. The graph of the approximate solution of the problem \eqref{Exam6} is depicted in Figure \ref{fig6.1}. \\
\section{Conclusion}
In this paper we have proposed a novel method for investigating the existence and uniqueness of a solution of a fully fourth order equation. It is based on the reduction of the original problem to an operator equation for the nonlinear term and the unknown values of the second derivatives $u''(0), u''(1)$. Under some easily verified conditions we have proved that the problem has a unique solution and it may be found by an iterative method. The convergence of the method as a geometric progression is established. At each iteration it is needed to solve two simple linear boundary value problems for second order equations. Several examples, where the exact solutions of the problem are known or  are not known, illustrated the effectiveness of the proposed method. Especially, some of the examples show the advantage of our results over ones of Agarwal. 

\par The proposed method can be applied to some other problems for ordinary and partial differential equations.

\section*{Appendices}
{\bf{A. 5-point numerical differentiation formula for the first derivative}} (see \cite{Li})
\begin{equation*}\label{eqI.1}
\left\{
\begin{array}{ll}
f'(x_0)=\dfrac{1}{12h}(-25f_0+48f_1-36f_2+16f_3-3f_4)+O(h^4), \\ 
f'(x_0)=\dfrac{1}{12h}(-3f_{-1}-10f_0+18f_1-6f_2+f_3)+O(h^4), \\ 
f'(x_0)=\dfrac{1}{12h}(f_{-2}-8f_{-1}+8f_1-f_2)+O(h^4), \\ 
f'(x_0)=\dfrac{1}{12h}(-f_{-3}+6f_{-2}-18f_{-1}+10f_0+3f_1)+O(h^4), \\ 
f'(x_0)=\dfrac{1}{12h}(3f_{-4}-16f_{-3}+36f_{-2}-48f_{-1}+25f_0)+O(h^4).
\end{array}
\right. 
\end{equation*}
{\bf{B. An approximation of order 4 on a solution of equation $u''=-f$}}  (see \cite{Sam}, Sect. 2.2)
\par In the one-dimensional, suppose that $u(x)$ is a solution of equation
\begin{equation} \label {eqI.2}
u''(x) = -f(x).
\end{equation}
Then the difference equation 
\begin{equation*}\label{eqI.3}
u_{\overline x x}=-f-\dfrac{h^2}{12}f_{\overline x x}
\end{equation*}
provides an approximation of fourth order on a solution $u(x)$ of equation \eqref{eqI.2}, where  
\begin{equation*}\label{eqI.4}
g_{\overline x x}=\dfrac{g(x+h)-2g(x)+g(x-h)}{h^2}.
\end{equation*}
{\bf{C. The Cavalieri-Simson Formular}} (see \cite[Sect. 9.2]{Quar})
Let $f$ be a real integrable function over the interval $[a,b]$. Introducing the quadrature nodes $x_k=a+kH/2$, for $k=0,1,...,2m$ and letting $h=\dfrac{b-a}{m}$ with $m\geq 1$. Then the composite Simpson formular is 
\begin{equation*}\label{eqI.5}
\begin{aligned}
&\int_a^b f(x)dx \approx \frac{h}{6} \Big [ f(x_0)+2\sum_{r=1}^{m-1}f(x_{2r})+4\sum_{s=0}^{m-1}f(x_{2s+1})+f(x_{2m})\Big]. 
\end{aligned}
\end{equation*}
The quadrature error associated with \eqref{eqI.5} is  $-\dfrac{b-a}{180}(h/2)^4f^{(4)}(\xi ), \quad \xi \in (a,b)$.

\end{document}